\documentclass[10pt, a4wide]{amsart}
\usepackage{graphicx}
\usepackage{amssymb,amsmath}
\usepackage{epsfig}
\usepackage{epsf}
\usepackage{float}
\usepackage{a4wide}
\usepackage{graphpap}
\usepackage{setspace}

\newtheorem{teo}{Theorem}[section]
\newtheorem{defi}[teo]{Definition}
\newtheorem{eje}[teo]{Example}
\newtheorem{lem}[teo]{Lemma}
\newtheorem{pro}[teo]{Proposition}
\newtheorem{coro}[teo]{Corolary}
\newtheorem{obs}[teo]{Remark}
\newtheorem{fact}[teo]{Fact}

\def\bteo{\begin{teo}}
\def\eteo{\end{teo}}
 \def\bdefi{\begin{defi}}
\def\edefi{\end{defi}}
\def\beje{\begin{eje}}
\def\eeje{\end{eje}}
\def\blem{\begin{lem}}
\def\elem{\end{lem}}
\def\bpro{\begin{pro}}
\def\epro{\end{pro}}
\def\bcoro{\begin{coro}}
\def\ecoro{\end{coro}}
\def\bob{\begin{obs}}
\def\eob{\end{obs}}
\def\bfact{\begin{fact}}
\def\efact{\end{fact}}

\def\E{{\mathbb E}}
\def\P{{\mathbb P}}
\def\R{{\mathbb R}}

\def\square{\ifmmode\sqr\else{$\sqr$}\fi}
\def\sqr{\vcenter{
         \hrule height.1mm
          \hbox{\vrule width.1mm height2.2mm\kern2.18mm\vrule width.1mm}
         \hrule height.1mm}}                  

\def\R{{\mathbb R}}  
\def\P{{\mathbb P}}  
\let\cal=\mathcal

\def\BB{{\cal B}}
\def\TT{{\cal T}}
\def\QQ{{\cal Q}}

\begin{document}
\title{A distance based test on random trees}
\author{Ana Georgina Flesia }
\address{Ana Georgina Flesia\\FaMAF-UNC\\Ing. Medina Allende s/n, Ciudad Universitaria\\ CP 5000, C\'ordoba, Argentina. }
\email{flesia@mate.uncor.edu}
\urladdr{http://www.famaf.unc.edu.ar/~flesia}
\thanks{ AGF was partially supported by PICT 2005-31659. \\Submitted to MACI2007, I CONGRESS ON COMPUTATIONAL, INDUSTRIAL AND APPLIED MATHEMATICS, October 2-5 2007. C\' ordoba, Argentina }
\author{Ricardo Fraiman}
\address{Ricardo Fraiman\\Departamento de Matem\'atica y Ciencias \\Universidad de San
  Andres\\ Buenos Aires, Argentina}
\email{rfraiman@udesa.edu.ar}


\begin{abstract}
In this paper, we address the question of comparison between populations of trees. We study an statistical test based on the
distance between empirical mean trees, as an analog of the two sample z statistic for comparing two means. Despite its simplicity, we can report that the test is  quite powerful to
separate distributions with different means but it does not
distinguish between different populations with the same mean, a more complicated test should be applied in that setting. The performance of the test is studied via
simulations on Galton-Watson branching processes. We also show an application  to a real data problem in genomics.
\end{abstract}
\maketitle
\section{Introduction}
\label{s1}

Random trees have long been an important modeling tool. Trees are useful when a
collection of observed objects are all descended from a common ancestral object
via a process of duplication followed by gradual differentiation. This would be the case of two
broad approaches to constructing random evolutionary trees: forwards in time
``branching process'' models, such as the Galton-Watson process, and
backwards-in-time ``coalescent'' models such as Kingman's coalescent (Kingman,
1982). We will show in our examples that the presence of specific short sequences or motifs in a string of elements taken from a finite alphabet is also related to a tree structure, so a random distribution of strings is directly related to a random distribution of trees.

In this preprint, we consider trees that have a root and evolve forward in time in discrete
generations, and each parent node (or vertex) having up to $m$ offspring nodes
in the next generation, as in Balding el al(2004), BFFS from now on. Given a suitable metric, BFFS prove law of large
numbers for empiric samples of trees and an invariance principle on the space of
continuous functions defined on the space of trees.

In this context, let  $\nu, \nu^*$ be distributions that give mass only to finite trees. The goal is
to test differences between the population laws
\begin{equation}
\label{test0}
 H_0: \nu = \nu^*\qquad H_A: \nu \neq \nu^*
\end{equation}
using i.i.d. random samples with distribution $\nu$ and $\nu^*$ respectively. Intuitively, if the expected mean of each population is different,
a naive test for this problem will reject the null hypothesis when the distance
between the empirical means associated with each sample is large enough, but it will fail if the population have different laws but the same expected mean. A Kolmogorov-type of test have been devised for this problem  in BFFS (2004) but a direct approach to calculate effectively the test statistic is quite difficult, since it is based on a supremo defined over
the space of all trees, which grows exponentially fast. The computation of the BFFS test, along with some discussion of identifiability of the measure have been worked out in Busch et al (2007).

In note we have studied the naive distance based test over simulations of Galton Watson processes, and we will also report an application to structural genomics, that is related to Variable Length Markov Chain Modeling. This is a similar example to the one introduced in Busch et al (2007), with another database, that relates to the work of Bejerano(2004).

\section{Trees, distances and tests}
\label{PyP}

We will review the definition of tree, that can be roughly thought as a set of nodes satisfying the condition
"son present implies father present". Let consider an alphabet $\cal A=\{1,\dots,m\}$, with $m\ge
2$ integer, representing the maximum number of children of a given node of the
tree. Let $V= \{1,2,\dots,m, 11,21,\dots,m1,12,22,32,\dots\}\cup \{\lambda\}$, the set of finite
sequences of elements in $\cal A$, plus the symbol $\lambda$ which
represents the {\em root of the tree}. The \emph{full tree} is the oriented
graph $ t_f=( V, E)$ with edges $ E\subset
 V\times V$ given by $ E=\{(v,av)\,:\,
v\in V,v\neq \lambda,\,a\in\cal A\}\cup \{(\lambda,a)\;a\in\cal A\}$, where $av$ is the sequence obtained by
concatenation of $v$ and $a$. In the full tree each node (vertex) has exactly
$m$ outgoing edges (to its offsprings) and one ingoing edge (from her father),
except for the root who has only outgoing edges.  The node $v=a_{k-1}\dots a_1$ is
said to belong to the \emph{generation} $k$; in this case we write ${\rm
  gen}(v)=k$. Generation 1 has only one node, the root.

We define a \emph{tree} as a function $t: V\to\{0,1\}$
satisfying
\begin{eqnarray}\label{i1}
t(v)&\ge&t(av).
\end{eqnarray}
for all $v\in V$ and $a\in A$, including the case of the root
\begin{eqnarray}\label{ii}
t(\lambda)&\ge&t(a).
\end{eqnarray} Abusing
notation, a tree $t$ is identified with the subgraph of the full tree $t=(V_t,E_t)$ with
\begin{eqnarray}
  \label{p1}
&& V_t=\{v\in V\,:\, t(v)=1\} \hbox{
and }
E_t = \{(v,av) \in E\,:\, t(v)= t(av)=1\}\;.
\end{eqnarray}
In figure~\ref{fig:vlmc1} we can observe a tree of depth 4. With this type of notation, the father of a node is written as a suffix in the description of the son, as it is often done in the definition of a Variable Length Markov Chain. We should notice though that the depth of a tree considers the root, and in VLMC models, the depth is the maximum length of a context, which do not consider the root.
\setlength{\unitlength}{1mm}
\begin{figure}[h]
  \begin{center}
    \begin{minipage}{5cm}
       \begin{picture}(50,68)(-2,-7)
        \put(1,-7){\includegraphics*[height=7cm]{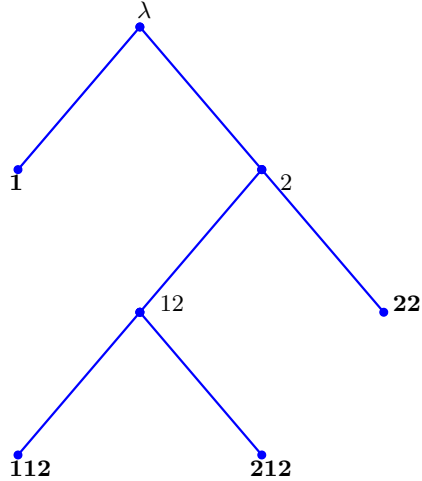}}
        \put(25,59){\small $\lambda$}
        \put(44,36){\small${2}$}
        \put(8,36){\small$\mathbf{1}$}
        \put(28,20){\small${12}$}
        \put(59,20){\small$\mathbf{22}$}
        \put(8,-2){\small$\mathbf{112}$}
        \put(40,-2){\small$\mathbf{212}$}
      \end{picture}
    \end{minipage}
     \caption{An example of tree of depth 4. The leaves are written in boldface. }
     \label{fig:vlmc1}
     \end{center}
   \end{figure}

 Let $\TT$ be the set of all trees, and let $\phi: V\to\R^+$ be a strictly positive function such that
$\sum_{v\in V} \phi(v)<\infty$. We define a distance between two trees in $\TT$ as a weighted sum over the nodes that are present in a tree and absent in the other, following the formula
\begin{eqnarray}
d(t,y)=\sum_{v\in V}\phi(v)|t(v)-y(v)|.
\end{eqnarray}
as it have been done in BFFS(2004).
The natural sigma algebra is the minimal one containing cylinders, sets of trees
defined by the presence/absence of a finite number of nodes. The natural
topology is the one generated by the cylinders as open sets. So it is easy to prove that the distance $d$ we defined before  generates the natural topology, and  $(\TT,d)$ becomes a compact metric space, see BFFS (2004).  We denote $\BB$
the $\sigma$-field of Borel subsets of $\TT$, induced by the metric $d$.

\paragraph*{\bf Random trees}
A \emph{random tree} with distribution $\nu $ is a measurable
function
\begin{equation}
  \label{p3}
  T:\Omega\to\TT\;\hbox{such that}\; \P(T\in A) = \int_A \nu(dt) \;.
\end{equation}
for any Borel set  $A\in \BB$, where $(\Omega,\mathcal F,\P)$ is a
probability space and $\nu$ a probability on $(\TT,\BB)$.

\paragraph*{\bf Expected value}

 The expected value or $d$-mean of a random tree $T$ is the set (of trees) $\E_d
 T$ which minimizes the expected distance to $T$:
 \begin{eqnarray}
  \label{p4}
  \E_d T := \arg \min_{t\in \TT}\, \int_\TT d(t,y) \;\nu(dy)
 \end{eqnarray}

The set $\E_d T$ is not empty, see BFFS (2004). Any element of the set $\E_dT $ is
also called a $d$-mean or $d$-center. Since $\E_dT$ depends only on the
distribution $\nu$ induced by $T$ on $\TT$, it may also be denoted as
$\E_d(\nu)$.

\paragraph*{\bf Empiric mean trees} Let ${\bf T}=(T_1,\dots,T_n)$ be a random
sample of $T$ (independent random trees with the same law as $T$). The empiric mean tree (empiric $d$-center,
 sample $d$-mean) is defined as the random set of trees given by
 \begin{equation}
   \label{p5}
   \overline {\bf T} := \arg\min_{t\in\TT} \frac1n\sum_{i=1}^nd(T_i,t)
 \end{equation}
 This formula may show the problem as more difficult that it is, since it is calling for a search over the whole set of trees, that grows exponentially in the number of nodes. But it is easy to prove that the empiric mean tree of a set of trees can be built by majority vote over the nodes. That means, at least one of them can be defined as the tree whose nodes are present only if they are present in at least half of the sample.
 \begin{pro}
 Let ${\bf T}=(T_1,\dots,T_n)$ be a random
sample of $T$ , and let $t^*$ be the tree defined as the tree whose nodes are present only if they are present in at least half of the sample. Then $t^*$ is an empiric mean tree.
\end{pro}
{\bf Proof}
 Let first notice that if $t\in \TT$

 \begin{eqnarray*}
 \frac1n\sum_{i=1}^nd(T_i,t)&=& \frac1n\sum_{i=1}^n\sum_{v\in V}\phi(v)|T_i(v)-t(v)|\\&=&\sum_{v\in V_t}\phi(v)\frac1n\sum_{i=1}^n|T_i(v)-t(v)|+\sum_{v\in \cup V_{T_i}/V_t}\phi(v)\frac1n\sum_{i=1}^n|T_i(v)-t(v)|\end{eqnarray*}

 \begin{eqnarray*}
 &=&\sum_{v\in V_t}\phi(v)\frac{\# \mbox{trees in the sample $v$ is not present}}{n}+
\sum_{v\in \cup V_{T_i}/V_t}\phi(v)\frac{\# \mbox{trees in the sample $v$ is present}}{n}\\
 \end{eqnarray*}

 So, to reduce the average of distances we have to reduce both summands, keeping and adding nodes to the candidates of empiric means. The first point to notice is that the first summand is reduced when the candidate $t$ keeps nodes that are present in many trees of the sample. If $t$ keeps a node that is not in any tree of the sample, the first summand adds the full value of $\phi(v)$. The second summand is reduced when the tree $t$ do not keep a node that is present only in a few trees of the sample. The cut off that balance the presence-absence relationship for each node is then $1/2$.\hfill\square

  \begin{obs}
  We should notice that if the number of trees in the sample is odd, the empiric mean is unique, but if the sample size is even, the node that it is present in exactly a half of the sample can be kept or not, without increasing the distance, so we will have at least two empiric means,  one will have the least of possible present nodes, and the other the most.
\end{obs}

\paragraph*{\bf Example 1: Galton-Watson related population of trees}

We consider binary trees, $m=2$, the extension to an arbitrary number of offsprings
$m$ is straightforward. In a binary binomial Galton-Watson model, the offspring
number is $0,1$ or $2$ with probabilities $(1-p)^2,\ 2p(1-p)$ and $p^2$.  The expected
mean tree keeps a node $v$, if and only if ${\rm gen}(v)\leq k_0$,
where $k_0 = \max\{k\in \{0,1\dots\}:p^k\geq 1/2\}$.  When $p<1/2$, the expected mean
tree is the empty tree. For instance, if $p=0.5$ and $p^* =
0.75$, the expected mean trees are $T_p =\{\lambda\}$ and
$T_{p^*}=\{\lambda,1,2\}$, the full trees of depth 1 and
2 respectively, but for $p\in [0.5,0.70]$ the population have the same expected mean tree. This is a very simple parametric case where the maximum likelihood test has maximum power, so it is not of much use to introduce a new test in this setting, {\em if we knew } that we have a Galton Watson process producing our observations. We consider this example only  to asses the power of the proposed test via simulation.

\paragraph*{\bf Example 2: Variable Length Markov Chains and related population of trees}

A Variable Length Markov Chain is a stochastic process introduced first  by Rissanen (1983) in the setting of information theory, and that have been recalled lately by B\"uhlmann and Wyner (1999), and many others in the context of Protein Functionality Modeling, see Bejerano (2003) and references therein.
 In this model the probability of occurrence of each symbol at a given time depends on a
finite number of precedent symbols.  The number of relevant precedent symbols
may be variable and depends on each specific sub-sequence.  More precisely, a
VLMC is a stochastic process $(X_n)_{n\in\mathbb Z}$, with values on a finite
alphabet $\cal A$, such that
  \begin{equation}\label{eq:prob}
    P[X_n=\cdot \,|\, X_{-\infty}^{n-1}=x_{-\infty}^{n-1}] = P[X_n=\cdot
    \,|\,X_{n-k}^{n-1}=x_{n-k}^{n-1}]\,,
\end{equation}
where $x_{s}^{r}$ represents the sequence $x_s,x_{s+1},\dotsc,x_r$ and $k$ is a
stopping time that depends on the sequence $x_{n-k},\ldots,x_{n-1}$.  As the
process is homogeneous, the relevant past sequences $(x_{n-k},\ldots,x_{n-1})$ do
not depend on $n$ and are  called \emph{contexts}, and denoted by $(x_{-k},\ldots,x_{-1})$. The set of all contexts $\tau$
can be represented as a rooted tree $t$, where each complete path from the
leaves to the root in $t$ represents a context. Calling $p$ the transition
probabilities associated to each context in $\tau$ given by (\ref{eq:prob}), the
pair $(\tau,p)$, called \emph{probabilistic context tree}, has all information
relevant to the model, see Rissanen (1983)  and B\"uhlmann and Wyner (1999).
As an example, take a binary alphabet $\cal A = \{1,2\}$ and transition probabilities
   \begin{equation}
\label{e27}
  P[X_n=x_n \,|\, X_{-\infty}^{n-1}=x_{-\infty}^{n-1}] =
  \begin{cases}
    P[X_n=1 \,|\, X_{n-2}^{n-1} = 1\ 1]& = 0.7, \\
    P[X_n=1 \,|\, X_{n-2}^{n-1} = 2\ 1 ]&=
      0.4, \\
    P[X_n=1 \,|\, X_{n-1}=2]&=0.2.
  \end{cases}
\end{equation}
so that, if $x_{n-1}=2$, then the stopping time $k=1$ and $X_n=1$ with probability $0.2$;
otherwise the stopping time is $k=2$ and $X_n=1$ with probability $0.7$ if both $x_{n-1}=x_{n-2}=1$
or with probability $0.4$ if $x_{n-1}=1$ and $x_{n-2}=2$. The set of contexts is
$\tau = \{11,21,2\}$, when the set of all active nodes of the associated rooted tree $t$ is $V_t = \{1,11,21,2,\lambda\}$, since $1$ is an internal node in the path of the context $11$ and $21$, and $\lambda$ is the root. Another example over the same alphabet is given by the transition probabilities
\begin{equation}
\label{e28}
  P[Y_n=y_n \,|\, Y_{-\infty}^{n-1}=y_{-\infty}^{n-1}] =
  \begin{cases}
    P[Y_n=1 \,|\, Y_{n-1} = 1]& = 0.6, \\
    P[Y_n=1 \,|\, Y_{n-2}^{n-1} = 2\ 2 ]&=
      0.4, \\
    P[Y_n=1 \,|\, Y_{n-2}^{n-1}=1\ 2]&=0.2.
  \end{cases}
\end{equation}
The set of contexts is
$\eta = \{1,12,22\}$, when the set of all active nodes of the rooted tree $y$ is $V_y = \{\lambda,1,12,2,22\}$, since $2$ is an internal node in the path of the context $12$ and $22$.
The corresponding rooted trees $t$ and $y$ are
represented in Figure~\ref{fig:vlmc}. Let compute the distance between the these two trees,
\begin{eqnarray*}
d(t,y)&=&\sum_{v\in V}\phi(v)|t(v)-y(v)|\\&=&\phi(\lambda)|t(\lambda)-y'\lambda)|
+\phi(1)|t(1)-y(1)|+\phi(21)|t(2)-y(2)|+\phi(11)|t(11)-y(11)|\\&&
+\phi(12)|t(12)-y(12)|+\phi(21)|t(21)-y(21)|+\phi(22)|t(22)-y(22)|\\
&=&0+0+0+\phi(11)
+\phi(12)+\phi(21)+\phi(22)\\
&=&4\times 0.36^3=0.186624
\end{eqnarray*}
considering $\phi(v)=z^{gen(v)}$, $z=0.36$.
\setlength{\unitlength}{1mm}
\begin{figure}[t]
    \begin{minipage}{3cm}
       \begin{picture}(50,58)(-2,-7)
        \put(-18,47){(a)}
        \put(-17,5){\includegraphics*[height=3.5cm]{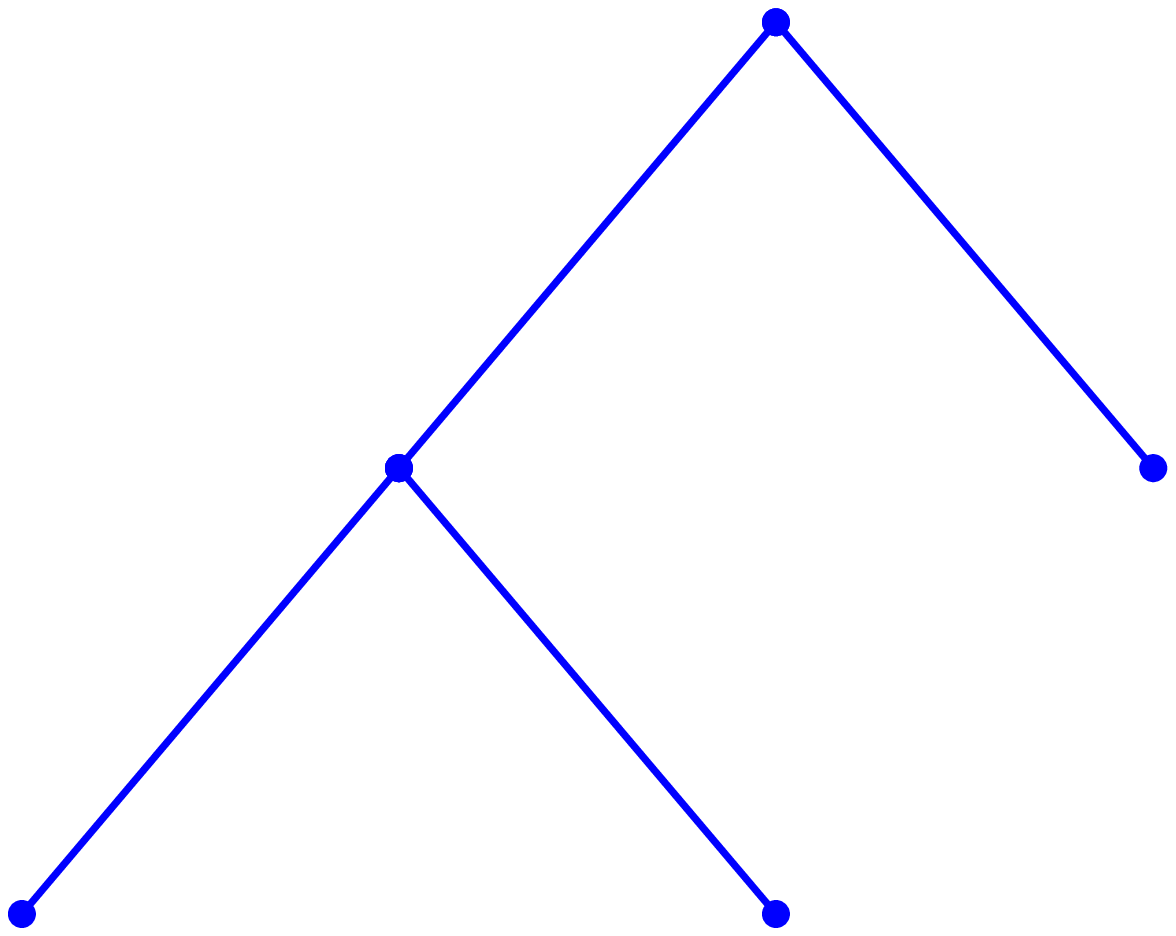}}
        \put(12,38.5){\small $\lambda$}
        \put(26,22){\small$\mathbf{2}$}
        \put(26,18){\scriptsize $(0.2,0.8)$}
        \put(-5,22){\small 1}
        \put(15,5){\small$\mathbf{21}$}
        \put(15,0){\scriptsize $(0.4,0.6)$}
        \put(-15,5){\small$\mathbf{11}$}
        \put(-16,0){\scriptsize $(0.7,0.3)$}
      \end{picture}
    \end{minipage}
    \hspace{3cm}
    \begin{minipage}{3cm}
      \begin{picture}(50,58)(-2,-7)
      \put(-18,47){(b)}
        \put(-17,5){\includegraphics*[height=3.5cm]{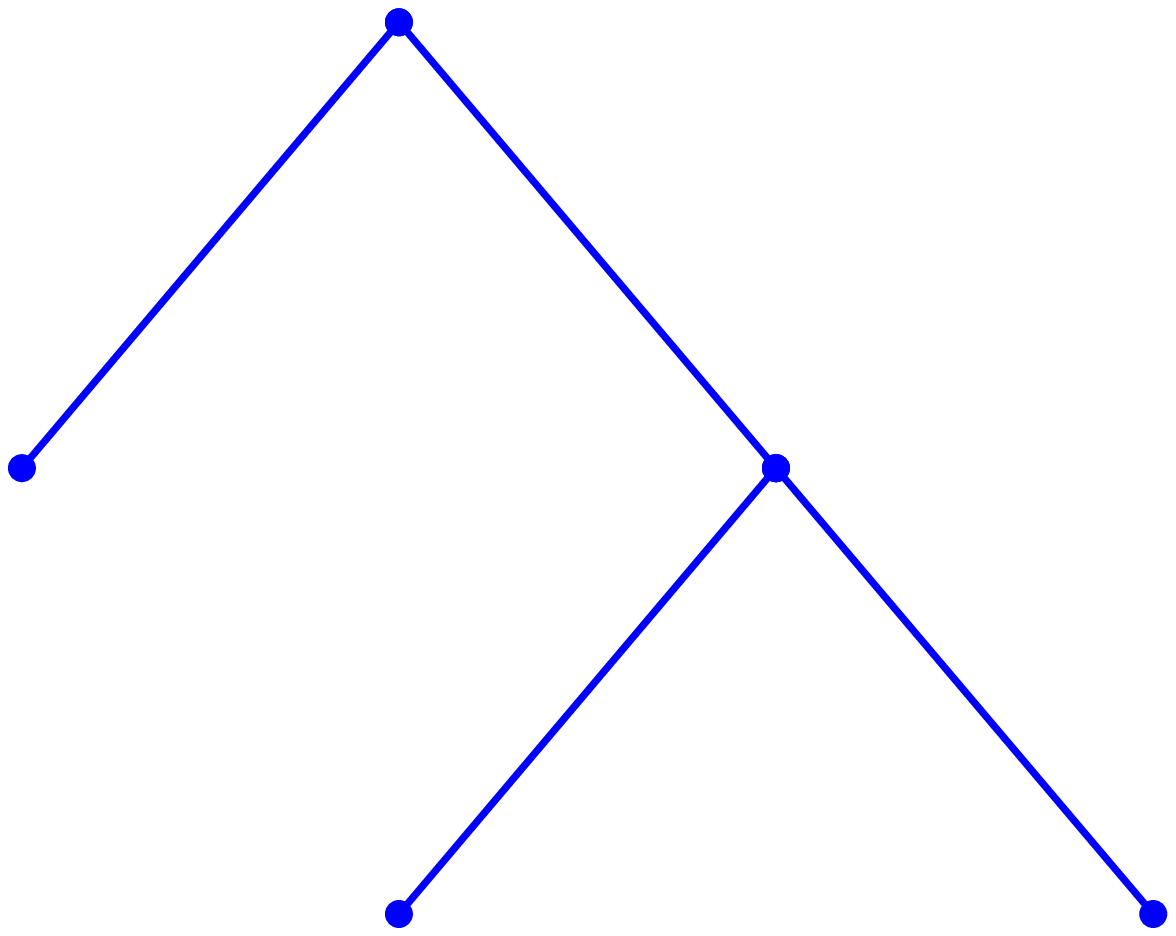}}
        \put(0,38.5){\small $\lambda$}
        \put(16,22){\small 2}
        \put(-15,18){\scriptsize $(0.6,0.4)$}
        \put(-15,22){\small $\mathbf{1}$}
        \put(-5,6){\small$\mathbf{12}$}
        \put(22,0){\scriptsize $(0.4,0.6)$}
        \put(26,6){\small$\mathbf{22}$}
        \put(-8,0){\scriptsize $(0.2,0.8)$}
      \end{picture}
    \end{minipage}
     \caption{An example of two probabilistic context trees over the alphabet
       $A=\{1,2\}$. (a) The tree $t$
       represents the pair $(\tau,p)$, where $\tau = \{11,21,2\}$ is the set of contexts and $p$ are
       the transition probabilities given by (\ref{e27}).(b) The tree $y$
       represents the pair $(\eta,q)$, where $\eta = \{12,22,1\}$ is the set of contexts and $q$ are
       the transition probabilities given by (\ref{e28}). }
     \label{fig:vlmc}
   \end{figure}

Now, let suppose that we are given a sequence of symbols that have been produced by a VLMC with an unknown context tree. There are several algorithms that estimates the context tree associated to the chain using the sequence as an input. Let fix the rule of estimation, par example, the Probabilistic Suffix Trees algorithm (PST) from Bejerano (2004). This rule is a {\em random tree} that generate trees in $\TT$ following a given probability distribution $\nu$ that is associated to the chain. If we have two independent samples of strings that have been hypothetically produced by two different unknown chains, we would like to derive a test that will rule if there is evidence in the samples to support that hypothesis. We should stress the fact that we are not using the probability transitions but the structure of the estimated context trees to derive the test. In the case that the chains have the same structure and the probability transitions are different, this approach will not apply.

Busch et al (2007) go further in this line of reasoning, suggesting a test that can rule when two samples from a collection of different VLMV models, clustered by a specific characteristic, are significantly different or not. Our test is not as general as it, but it is very simple to understand and to compute, and this ideas can extend easily to clustering and discrimination problems, that are based on distance. In Flesia et al (2007), we are currently working with an extension of the K-means algorithm for clustering, and k nearest neighbors procedure for discrimination, with a population of trees that are estimations of a VLMC context tree chain.

\paragraph*{\bf Testing differences of populations}
\label{testing} We consider measures $\nu\in\QQ_f$, the space of
probability measures that concentrate mass on trees with a finite number of
nodes. We describe the two-sample problem.

Let $\nu, \nu^*$ be distributions in $\QQ_f$. The goal is
to test
\begin{equation}
\label{test00}
 H_0: \nu = \nu^*\qquad H_A: \nu \neq \nu^*
\end{equation}
using i.i.d. random samples ${\bf T}=(T_1,\dots,T_n) $ and ${\bf
  T}^*=(T_1^*,\dots,T_m^*)$ with distribution $\nu$ and $\nu^*$ respectively.

\paragraph*{\bf Test based on the distance between mean trees}

When the expected $d$-means are different, $\E T\neq \E T^*$, one expects that the
distance between the empirical mean trees $\overline {\bf T}, \overline {\bf
  T}^*$ will be positive, for functions $\phi$ which do not penalize too much
the first generations, as $\phi(v)=z^{{\rm gen}(v)}$ with $0<z<\frac1m$. A simple
and naive test for this problem will reject the null hypothesis when the distance
between the empirical means associated with each sample is large enough.

\paragraph*{\bf Computation}
The lack of knowledge of the distribution of the distance between empirical
means may be overcame using Monte Carlo randomization.  If the null hypothesis is
$\nu=\nu^*$,  and
  \[d_k= d(\overline {\mathbf T},\overline {\mathbf T}^*)
= \sum_{v\in V} |\overline {\mathbf T}(v)-\overline {\mathbf T}^*(v)| \phi(v)
\]
 is the empiric distance between the mean trees of the $k$th pair of
simulated sample, created by randomly rearranging the whole set of observations, and assigning the first $n_1$ observations to the first sample and the rest to the second sample, we define the quantile $q_{\alpha}$ as the value such that
\[
\alpha = P(d(\overline{{\mathbf T}},\overline{\mathbf T}^* )>
q_{\alpha})
\]
This value can be approximated using the order statistics $ d^{(1)},\dots,d^{(N)}$ and
taking $q_{\alpha}$ as $d^{([N(1-\alpha)])}$ (here $[a]$ denotes the greatest
integer not greater than $a$).  For the original samples ${\mathbf T}$ and ${\mathbf
  T}^*$, the test will reject the hypothesis if $d(\overline{\mathbf
  T},\overline{\mathbf T}^* )>q_{\alpha}$ at level~$\alpha$.  The type-2 error
can be estimated analogously for each alternative hypothesis $\nu_a$.

\section{Computational examples}

\paragraph*{\bf Simulation}
 To study the
performance of the tests on a controlled environment we simulate several
populations of trees using Galton-Watson processes and simple variations of it.
We carefully choose the parameters to challenge the power of the tests.

Assume we have two random samples, each one from a Galton-Watson process with
possibly different parameters $p$ and $p^*$, denoted $GP(p)$ and
$GP(p^*)$. We we would like to test if these samples come from the same
process, that is,
\[
H_0: T\sim GP(p),\;  T^*\sim GP(p^*)\quad p=p^*\qquad
 H_A: T\sim GP(p),\; T^*\sim GP(p^*),\quad p\neq p^*
\]

In our simulation we already know the parameters of the underlying  distributions $\nu$ and $\nu^*$.  Thus, we have performed a Monte Carlo simulation test sampling trees from a mixture of both laws at random, until we reach the size of the first sample and label it
  sample 1. Then continue selecting with the same mixture, until we reach the
  size of the second sample, and label it sample~2. We compute the test statistics with these random samples, and store it, and repeat the process 1000 times. Then we generate  a fixed number of times a  sample from the distribution $\nu$, and a
sample from the distribution $\nu^*$, and calculate the test
statistics with them. If the true test statistic is greater than (1-$\alpha$)\% of the random values, then  the null hypothesis is rejected at $p<\alpha$.
The percentage of rejections for each value of $\alpha$ is considered  a measure of the power of the test.

We have computed the percentage of
rejection over 1000 tests of level $\alpha= 0.10, 0.05, 0.01$,
when $T^*_{k,1},\dots,T^*_{k,n}$ is $GP(p^*)$, with $p^* =0.6,
0.75, 0.8 $ and $0.85$, for sample sizes $n=31,51,101, 151$ and $201$.  The results are
reported on Table \ref{Tabla1}.
\begin{table}[h]
  \centering
  \begin{tabular}{|c|c c c c c|}
    \hline
             $\alpha=0.1$&{$n=31$}&{$n=51$}&{$n=101$}&{$n=151$}&{$n=201$}\\
             \hline
  $p=0.6$&  5.6 &    2.1  &        0  &        0  &        0 \\
  $p=0.75$& 52.8 &   65  &  92.5   & 99.3  & 100\\
  $p=0.8$&  86.2 &    93.7  &   99.9  &   100  &   100  \\
   $p=0.85$&99  &  99.9 &    100   & 100 &   100\\
  \hline
         $\alpha=0.05$&{$n=31$}&{$n=51$}&{$n=101$}&{$n=151$}&{$n=201$}\\
         \hline
   $p=0.6$&  5.60  &   02.1  &        0   &       0   &       0\\
   $p=0.75$& 52.80  & 47  & 92.5   & 99.3   & 100\\
   $p=0.8$&  78.60  &   93.7  &   94.7   &  100  &  100 \\
    $p=0.85$&97.80  &  99.8  &  100   & 100  &  100\\
  \hline
  $\alpha=0.01$&{$n=31$}&{$n=51$}&{$n=101$}&{$n=151$}&{$n=201$}\\
   \hline
     $p=0.6$  &0.70  &  2.10   &      0   &      0 &        0 \\
    $p=0.75$& 39.10   & 47.00    &58.40   & 95.10   & 96.9\\
     $p=0.8$  &51.70  &  76.90   & 94.70   & 95.10 &   96.2 \\
     $p=0.85$& 55.40  &  98.40  &  100  &  100 &   100\\
      \hline
  \end{tabular}
 \vspace{10mm} \caption{Percentage of rejections over 1000 tests, computed with with $p=0.5$ and $p^*=0.6,0.75,0.8,0.85$, sample size $n=31,51,101, 151$ and 201.}\label{Tabla1}
\end{table}

These results are in agreement with our intuitive ideas. As the sample size increases, the test is not able to reject the hypothesis of equal populations when $p=0.5$ and $p^*=0.6$, since their expected mean trees are equal. But when the expected mean trees are different, the test detects the difference with higher power as the sample size increases.

\paragraph*{\bf Variable Length Markov Chain Modeling of Protein Functionality}

A central problem in computational biology is to determine the function of a new discovered
protein using the information contained in its amino acid
sequence. Proteins are complex molecules composed by small blocks called amino acids. The amino acids
are linearly linked, forming a specific sequence for each protein. There exist 20 different amino acids represented by a one-letter code.

There are several problems related to protein functionality, we will only point out two of them here. One is the classification of the function of a new protein with the help of a training set, and the other is clustering a group of new and known proteins into meaningful functionality families. The goal of clustering protein sequences is to get a biologically meaningful partitioning. Genome projects are generating enormous amounts of sequence data that need to be effectively analyzed. Given to the amount of available data, and the lack of proper definition, clustering is a very difficult task,  so there is a need for ways of checking the validity of the partition proposed. As most databases are created by sequence alignment related methods, an impartial way of checking validity would be to apply an alignment free,  model based methodology.

 Most methods for clustering and classification need as input a similarity matrix, usually computed by sequence alignment. Model based clustering and classification without sequence alignment is leaded by Markov Chain modeling. Par example, Bejerano et al (2001) models protein sequences with stationary
Variable Length Markov Chains (VLMC), in order to classify a new given protein as belonging to the family whose model has higher probability of having produced that string. This approach needs also a reliable training set in order to build an accurate estimate of the unknown context tree of the chain.

 In this paper, we propose to check the coherence of of selected protein families performing a simultaneous hypothesis test, as it has been done in Busch et al(2007). We would like to test if several families that are members of a well known database are simultaneously significantly different. We are going to use the same database that was cited in Bejerano et al (2001), which provided the training data for the classification problem. The Pfam database is known to be a good reference for  protein functionality  clustering, so it would provide a benchmark for assessing the performance of our approach.

 We start modeling each functionality family of proteins as realizations of an unknown VLMC. But instead of learning the model using all the sequences of a given family to estimate the context tree with the Probabilistic Suffix Trees algorithm (PST) as in Bejerano et al (2001), we consider this rule as a {\em random tree} that generate one tree in $\TT$ per sequence. The  probability distribution $\nu$ of the random tree is associated to the chain that rules the family in an unknown fashion. If we have two independent samples of strings that have been hypothetically produced by two different unknown chains, we estimates with each of them the context tree of its chain and then consider we have two independent samples of trees, each one following a distribution associated to the family. We then test if there are enough evidence in the samples to reject the hypothesis of equal distribution. If we do reject, we consider the two families significantly different.

We must emphasize the difference between the example from Busch et al (2007) and our approach. They use the latest version of the Pfam database, which is significantly different from the one we are working, and they model each family as a collection of VLMC models, in correspondence with the notion of subfamily. We use the approach of Bejerano, modeling small families with only one VLMC, but estimating it several times with independent strings.

Let $\TT_4$ be the space of trees with $m=20$ possible children per node (the
symbols of the amino acid alphabet), and fixed maximum length $M=4$ and the parameter of the distance fixed as $z=0.36$. We test if ten
families selected from de P-fam database version 1, Bateman et al (2004), are simultaneously
significantly different using the following two step procedure
\begin{enumerate}
\item Transform the amino acid chains into trees via the Probabilistic Suffix
  Trees (PST) from Bejerano(2004), obtaining 10 samples of trees of
  maximum length of context equal to 3. The parameters of the PST have been set as the default.

\item Apply a Bonferroni correction to the 45 pairwise BFFS based comparisons,
  that means, each test is performed with a level of significance of
  $\alpha=0.05/45=0.001$ to get a simultaneous comparison of the 10 families,
  with overall level $\alpha=0.05$.
\end{enumerate}

We run all the pairwise tests at level $0.001$. We also run the tests
under the null hypothesis splitting each  data set at random in two subsets.
Table~\ref{t4} shows
the critical and the observed values for all pairwise tests of different
families (non--diagonal terms).  For the null hypotheses the observed value and
the p--value appear in boldface at the diagonal.  Despite the crude nature of
the Bonferroni method, the hypothesis of equal distribution is rejected in all cases when the samples came from different populations,
confirming the coherence of the selected protein families. In the case of the same family split in halves, we can observe p-values ranging from 0.12 to 0.90, values that can be used also to analyze the coherence of the family.

\begin{table}[h t]
\small
\label{t4}
\begin{center}
  \begin{tabular}{cccccc}
      \hline
      \hline
      Family &actin       &adh-short    &adh-zinc   &ank   &  ATP-synt-A    \\
      \hline
      actin  &{\bf(0.49,  0.71)}   &( 2.41,    9.85)&(3.37,  10.44)&(3.47,  9.84)&(4.62, 11.21)\\
      adh-short &  &{\bf(1.36,    0.43)} &(1.91, 4.91)&( 2.44,    5.55)&(2.05, 5.27)\\
      adh-zinc   & & &{\bf(1.66 ,   0.57)}&(2.58,   5.30)&(3.05,  6.70)\\
      ank&&&&{\bf(1.86  ,  0.81)}&(4.27,  8.37)\\
      ATP-synt-A &&&&&{\bf (1.67 0.52 )}\\\hline\hline
       Family     &beta-lactamase  &cox2 &cpn10 &DNA-pol &efhand\\
  \hline
  actin     & ( 3.76,  9.71)&(4.04,  11.46)&(5.52,   12.20)&(4.01,  9.73)&(3.06, 11.79)\\
  adh-short &(2.52 ,   3.91)&(2.14  ,  6.24)&(2.42 ,   6.20)&( 3.61 ,   6.44)&( 1.86,  6.07)\\
  adh-zinc  &(2.64  ,  5.14 )&(2.51 ,  7.38)&(2.79  , 7.91)&(2.58 ,  5.90)&(2.23, 6.93)\\
  ank       &(3.14,  5.32)&(3.03,  7.94)&( 5.04, 10.32)&(2.51, 3.16)&( 2.74, 8.79)\\
  ATP-synt-A   &(3.34,  6.25)&    (2.75,  4.95)&    (2.61,  5.28)&   ( 4.88,   8.98)&   (2.08,  6.27)\\
  beta-lactamase    &{\bf(1.819,   0.93)}&(2.98,   6.94)&  (  2.83,   6.52)&    (3.16, 6.58)&    (2.74,  6.49)\\
  cox2      && {\bf(2.05 ,   0.09)}&(2.95,  6.99 )&   (3.77,  9.19)&( 1.67,  6.49)\\
  cpn10     &  & &{\bf(1.30  ,  0.02)}&(6.46, 11.86)&( 2.02,   3.86)\\
  DNA-pol   &&&&{\bf(1.81 ,   0.23)}&(3.58 , 10.24)\\
  efhand    &&&&&{\bf( 0.65 ,   0.93)}\\\hline\hline
\end{tabular}
\end{center}
\caption{Critical value and observed value of 45 pairwise comparisons at level
  $\alpha=0.001$. Test rejects when the observed value is greater than the
  critical value. In boldface, observed value and  p-value when testing the same population, $N=1000$. The distance's  parameter zeta is equal to 0.35.}
\end{table}

\section{Final Remarks}

We have proposed a naive test to compare two population of trees with laws that do not have the same expected mean. The procedure is very simple, since it is based in the idea that the empiric mean tree of each sample, a strong consistent estimator of the expectation of the law that generates each population, should be separated in terms of BFFS distance. The test will reject the hypothesis of equal populations if the distance between the empiric means is big enough to ensure a small type one error. The quantile of the distribution has been derived by Monte Carlo randomization, and the power has been studied through Galton Watson simulations. We have also addressed a problem of functional genomics, to check the coherence of hypothesized functionality families. We suppose that each family of proteins is related to a random tree, and the allegedly members of each family form a sample of the law of the random tree that characterizes the family. We check if there is enough information in the samples to reject the hypothesis of  equal populations.

This approach will not work if the two populations have the same expected mean tree, as in the case of two sample of strings that have been generated by  chains with the same context tree but different transition probabilities. A more sophisticated test, the BFFS test, has already been proposed by  Balding et al (2004), a Kolmogorov type of test that maximizes the differences between the information of the samples, but it does not have a naive computation, since it involves a search over the set of trees that grows exponentially fast. In Busch et al (2007) the computation of the test has been derived and the performance of the test reported. Also, they suggest a way to model more complex group of proteins as collections of VLMC models, and test the same hypothesis with great success.

The key features of this test are the simplicity of the definition and it fast computation, that allows to realize easy preliminary approaches to the two samples testing problem.

\paragraph*{\bf Acknowledgments}
 I would like to thank Florencia Leonardi for providing the data used in our example of determination of protein functionality, which was also analized  in Leonardi (2007).

\end{document}